\documentclass[12pt,a4paper]{article}
\usepackage[latin2]{inputenc}
\usepackage{amsmath}
\usepackage{amsfonts}
\usepackage{amssymb}
\usepackage{graphicx}
\usepackage[left=2.00cm, right=2.00cm, top=2.50cm, bottom=2.50cm]{geometry}
\usepackage[T1]{fontenc}

\def\Q{{\mathbb Q}}
\def\Z{{\mathbb Z}}

\title{
Calculating generators of power integral bases in sextic fields with a quadratic subfield:
the general case
}
\author{
Istv\'an Ga\'al\\
{\small University of Debrecen, Mathematical Institute} \\
{\small H--4002 Debrecen Pf.400., Hungary,} \\
{\small e--mail: gaal.istvan@unideb.hu},
}

\begin{document}
\baselineskip=17pt

\maketitle
\thispagestyle{empty}

\renewcommand{\thefootnote}{\arabic{footnote}}
\setcounter{footnote}{0}

\vspace{0.5cm}


\begin{abstract}
In some previous works we gave algorithms for determining
generators of power integral basis in sextic fields with a quadratic subfield,
under certain restrictions.
The purpose of the present paper is to extend those methods to the general case,
when the relative integral basis of the sextic field over the quadratic subfield is of
general form. This raises several technical difficulties, that we consider here.
\end{abstract}

\section{Introduction}

Monogenity and power integral bases is a classical topic in algebraic number theory,
which is intensively studied even today, see \cite{nark} for classical results,
\cite{book} and \cite{axioms} for more recent results.

A number field $K$ of degree $n$ with ring of integers $\Z_K$ 
is called {\it monogenic} (cf. \cite{book}) if there 
exists $\xi\in \Z_K$ such that $(1,\xi,\ldots,\xi^{n-1})$ is an integral basis, 
called {\it power integral basis}. 
We call $\xi$ the {\it generator} of this power integral basis.

An irreducible polynomial $f(x)\in\Z[x]$ is called {\it monogenic}, 
if a root $\xi$ of $f(x)$ generates a power integral basis in $K=\Q(\xi)$.
If $f(x)$ is monogenic, then $K$ is also monogenic, but the converse is
not true.

For $\alpha\in\Z_K$ (generating $K$ over $\Q$)  the module index 
\[
I(\alpha)=(\Z_K:\Z[\alpha])
\]
is called the {\it index} of $\alpha$. The element 
$\alpha$ generates a power integral basis in $K$
if and only if $I(\alpha)=1$. 
Searching for elements of $\Z_K$, generating power integral bases, leads to
a Diophantine equation, called {\it index form equation} (cf. \cite{book}). 

There are certain  algorithms to determine "all solutions" of these equations, that is 
all generators of power integral bases. This "complete resolution" 
requires very often too long CPU time. On the other hand 
there are some very fast methods for 
determining generators of power integral bases with "small" coefficients,
say, being $<10^{100}$ in absolute value, with respect to an integral basis.
All our experiences show that generators of power integral bases have
very small coefficients in the integral basis, therefore 
these "small" solutions cover all solutions with high probability, 
certainly all generators that
can be used in practice for further calculations. 
It is usual to apply such algorithms also if we need to solve a large number of 
equation (cf. \cite{book}).

In sextic fields with a quadratic subfield we developed some efficient methods
for calculating "small" solutions of the index form equation, see 
\cite{grs62}, \cite{g25}. For simplicity, in these results we assumed, that
the basis of the sextic field is of special type
\[
(1,\alpha,\alpha^2,\omega,\omega\alpha,\omega\alpha^2),
\]
where $(1,\omega)$ is an integral basis of the quadratic subfield.
This implicitly yields, that the sextic field apriori has
a relative power integral basis over the quadratic subfield. 
In the present paper we extend this special case to the general case, when
the relative integral basis is of arbitrary form.

This paper was initiated by the recent work of 
Harrington and Jones \cite{hj6}, where they consider sextic trinomials of the form
$f(x)=x^6+ax^3+b$. Considering the sextic fields in \cite{hj6}, generated by a root of 
such a trinomial, we find that in most cases the root of the polynomial 
does not generate a power integral basis over the quadratic subfield.
In the present paper we intend to give a fast algorithm 
to calculate "small" solutions of the index form equation
in such sextic fields. We shall see, that 
some crucial ingredients of the method are similar to the formerly 
considered simpler cases, however several complications occur that make it worthy to provide a description in the general case. 
In other words, we describe how the previous algorithms 
can be extended to the general case. Also, note that the present method 
can be easily transformed to a process to calculate all solutions.
 
\section{Sextic fields with a quadratic subfield}

Let $M$ be a quadratic number field with integral basis $(1,\omega)$,
and let $f(x)=x^3+C_2x^2+C_1x+C_0\in\Z_M[x]$ be the relative defining
polynomial of $\alpha$ over $M$, with $K=M(\alpha)$. 
For sextic fields with a quadratic subfield a crucial step,
the reduction, only works for complex quadratic subfields, therefore
we assume that $M$ is complex. 

We are going to determine generators of power integral bases of $K$.

To present our formulas explicitly we write the relative integral basis 
of $K$ over $M$ in the form
\begin{equation}
\left(1,\frac{A\alpha+B}{k},\frac{C\alpha^2+D\alpha+E}{\ell}\right),
\label{basis}
\end{equation}
where $A,B,C,D,E\in\Z_M$, $0<k,\ell\in\Z$. 
Note that if $K$ is (absolute) monogenic, then it is also relative monogenic over $M$,
implying that $K$ has a relative integer basis over $M$.

Using the relative integral basis (\ref{basis}) we can represent any $\gamma\in\Z_K$
in the form
\begin{equation}
\gamma=X_0+X_1\frac{A\alpha+B}{k}+X_2\frac{C\alpha^2+D\alpha+E}{\ell},
\label{gamma}
\end{equation}
with unknown $X_i=x_{i1}+\omega x_{i2}\in\Z_M\; (i=0,1,2)$. 
Our purpose is to construct a fast algorithm 
to determine all tuples $(x_{02},x_{11},x_{12},x_{21},x_{22})\in\Z^5$ with
\begin{equation}
\max(|x_{02}|,|x_{11}|,|x_{12}|,|x_{21}|,|x_{22}|)<C,
\label{xkorl}
\end{equation}
with say, $C=10^{100}$,
such that $\gamma$ generates a power integral basis in $K$ 
(the index of $\gamma$ is independent from $x_{01}$).

We have
\[
\gamma=Y_0+Y_1\alpha+Y_2\alpha^2,
\]
where
\begin{equation}
Y_0=X_0+X_1\frac{B}{k}+X_2\frac{E}{\ell},\; Y_1=X_1\frac{A}{k}+X_2\frac{D}{\ell},\;
Y_2=X_2\frac{C}{\ell},
\label{Y}
\end{equation}
are not necessarily integer elements in $M$.

Let $\mu^{(i)}$ be the conjugates of any $\mu\in M$,
corresponding to $\omega^{(i)},\; (i=1,2)$.
We denote by $\alpha^{(i,j)}$ the roots of 
$f^{(i)}(x)=x^3+C_2^{(i)}x^2+C_1^{(i)}x+C_0^{(i)},\; (i=1,2,j=1,2,3)$.
The conjugates of any $\mu\in K$ corresponding to $\alpha^{(i,j)}$
will also be denoted by $\mu^{(i,j)}$.

For $i=1,2,1\le j_1,j_2\le 3, j_1\ne j_2$ we have
\[
\frac{\gamma^{(i,j_1)}-\gamma^{(i,j_2)}}{\alpha^{(i,j_1)}-\alpha^{(i,j_2)}}
=
Y_1+(\alpha^{(i,j_1)}+\alpha^{(i,j_2)})Y_2
=
Y_1-\delta^{(i,j_3)}\; Y_2,
\]
where $-\delta^{(i,j_3)}=\alpha^{(i,j_1)}+\alpha^{(i,j_2)}=-C_2^{(i)}-\alpha^{(i,j_3)}$,
$C_2\in\Z_M$ being the quadratic coefficient of the relative defining
polynomial $f(x)$ of $\alpha$ over $M$, and $\{j_3\}=\{1,2,3\}\setminus \{j_1,j_2\}$.

By the representation (\ref{basis}) of the relative integral basis of $K$
over $M$, for the relative discriminant $D_{K/M}$ we have 
\[
N_{M/\Q}(D_{K/M})=\frac{N_{M/\Q}(d_{K/M}(\alpha))}{(k\ell)^2}
=\frac{1}{(k\ell)^2}
\prod_{i=1}^2\prod_{1\le j_1<j_2\le 3}\left(\alpha^{(i,j_1)}-\alpha^{(i,j_2)}\right)^2.
\]

Therefore we obtain
\[
I_{K/M}(\gamma)=\frac{1}{\sqrt{|N_{M/\Q}(D_{K/M})|}}
\prod_{i=1}^2\prod_{1\le j_1<j_2\le 3}\left|\gamma^{(i,j_1)}-\gamma^{(i,j_2)}\right|
\]
\begin{equation}
=(k\ell)\prod_{i=1}^2\prod_{j=1}^3\left|Y_1^{(i)}-\delta^{(i,j)}Y_2^{(i)}\right|
=
(k\ell) \; |N_{M/\Q}(N_{K/M}(Y_1-\delta Y_2))|.
\label{thue}
\end{equation}

As it is known (see \cite{book}, Chapter 1, Theorem 1.6), if $I(\gamma)=1$, then both
\[
I_{K/M}(\gamma)=1,\;\; \text{and}\;\; J(\gamma)=1,
\]
where
\begin{equation}
J(\gamma)=\frac{1}{(\sqrt{|D_M|})^3}\prod_{j_1=1}^3\prod_{j_2=1}^3
\left|\gamma^{(1,j_1)}-\gamma^{(2,j_2)}\right|.
\label{JG}
\end{equation}

\section{Reduction}

By $I_{K/M}(\gamma)=1$, (\ref{thue}) implies 
\begin{equation}
N_{M/\Q}(N_{K/M}(Z_1-\delta Z_2))=\pm (k\ell)^5,
\label{T2}
\end{equation}
where
\begin{equation}
Z_1=(k\ell)Y_1,\;\; Z_2=(k\ell)Y_2\in\Z_M.
\label{ZY}
\end{equation}
Using an algebraic number theory package like 
Magma or Kash we can determine a complete set of non-associated
elements $\mu\in\Z_M$ of norm $\pm (k\ell)^5$. Let $\varepsilon$ 
be one of the finitely many units in $M$. We confer
\begin{equation}
N_{K/M}(Z_1-\delta Z_2)=\varepsilon\mu,
\label{mu}
\end{equation}
with certain possible values of $\mu,\varepsilon$.
In complex quadratic fields the conjugated elements have equal absolute values,
therefore (\ref{mu}) implies
\begin{equation}
\prod_{j=1}^3\left|Z_1^{(1)}-\delta^{(1,j)}Z_2^{(1)}\right|=|k\ell|^{5/2}.
\label{52}
\end{equation}
Denote by $j_0$ the conjugate with
\[
\left|Z_1^{(1)}-\delta^{(1,j_0)}Z_2^{(1)}\right|
=\min_{1\le j\le 3}\left|Z_1^{(1)}-\delta^{(1,j)}Z_2^{(1)}\right|.
\]
Then
\begin{equation}
\left|Z_1^{(1)}-\delta^{(1,j_0)}Z_2^{(1)}\right|\le c_1
\label{min}
\end{equation}
with $c_1=|k\ell|^{5/6}$,
and for $j\ne j_0$ we have
\begin{equation}
\left|Z_1^{(1)}-\delta^{(1,j)}Z_2^{(1)}\right|\ge 
\left|\delta^{(1,j)}-\delta^{(1,j_0)}\right|\; |Z_2^{(1)}|-c_1
\ge c_2\;|Z_2^{(1)}|,
\label{nagy}
\end{equation}
with $c_2=0.9\cdot \min_{j\ne j_0}|\delta^{(1,j)}-\delta^{(1,j_0)}|$, if 
$|Z_2^{(1)}|>10c_1/\min_{j\ne j_0}|\delta^{(1,j)}-\delta^{(1,j_0)}|$.
Small coordinates of $Z_2$, not satisfying this inequality
are tested separately. 

We set $Z_1=z_{11}+\omega z_{12},Z_2=z_{21}+\omega z_{22}$ with 
$z_{11},z_{12},z_{21},z_{22}\in\Z$ and let
\[
A=\max(|z_{11}|,|z_{12}|,|z_{21}|,|z_{22}|).
\]
Note that to find all suitable $(x_{02},x_{11},x_{12},x_{21},x_{22})\in\Z^5$
satisfying (\ref{xkorl}), in view of (\ref{Y}), (\ref{ZY}) we have to consider 
all $(z_{11},z_{12},z_{21},z_{22})$ with 
\begin{equation}
A\le 6 (k\ell) C (1+\overline{|\omega|})\max\left(1,\frac{B}{k},\frac{E}{\ell},
\frac{A}{k},\frac{D}{\ell},\frac{C}{\ell}\right).
\label{aaa}
\end{equation}
(\ref{min}) implies 
\[
\max(|z_{11}|,|z_{12}|)\le 2|Z_1^{(1)}|
\le 2 (c_1+\overline{|\delta|} \;|Z_2^{(1)}|)
\le 2 (0.1+\overline{|\delta|}) \;|Z_2^{(1)}|,
\]
if $0.1 |Z_2^{(1)}|>c_1$ (small coordinates of $Z_2$
are tested separately). Here $\overline{|\delta|}$ is the size 
$\delta$ (the maximum absolute values of its conjugates).
Similary $\max(|z_{21}|,|z_{22}|)\le 2|Z_2^{(1)}|$, therefore
\[
A\le 2 (0.1+\overline{|\delta|}) \;|Z_2^{(1)}|.
\]

By (\ref{52}) and (\ref{nagy}) we obtain 
\[
\left|Z_1^{(1)}-\delta^{(1,j_0)}Z_2^{(1)}\right|
\le \frac{(k\ell)^{5/2}}{c_2^2} |Z_2^{(1)}|^{-2}\le c_4 A^{-2},
\]
with
\[
c_4=\frac{(k\ell)^{5/2}}{4 c_2^2 (0.1+\overline{|\delta|})^2 },
\]
whence
\begin{equation}
\left| z_{11}+\omega^{(1)}z_{12}
- \delta^{(1,j_0)}z_{21}-\delta^{(1,j_0)}\omega^{(1)}z_{22}\right|\le c_4 A^{-2}.
\label{forma}
\end{equation}

The bound in (\ref{aaa}) is reduced in several consecutive steps. We start with 
$A_0=A_{\max}$, $A_{\max}$ being the bound in (\ref{aaa}). We assign a suitable
large constant $H$, perform the following reduction step, which produces a new bound
for $A$. We set this new bound in place of $A_0$ and continue the reduction
until the reduced bound is smaller then the original one.

Consider the lattice generated by the
 columns of the matrix
\[
\left(
\begin{array}{cccc}
1&0&0&0\\
0&1&0&0\\
0&0&1&0\\
0&0&0&1\\
   H&H\Re(\omega^{(1)})&\Re(-\delta^{(1,j_0)})&H\Re(-\delta^{(1,j_0)}\omega^{(1)})\\
   0&H\Im(\omega^{(1)})&H\Im(-\delta^{(1,j_0)})&H\Im(-\delta^{(1,j_0)}\omega^{(1)})\\
\end{array}
\right).
\]
Denote by $b_1$ the first vector of the LLL reduced basis of this lattice.
According to Lemma 5.3 of \cite{book}, if $A\le A_0$ and $H$ is large enough to have
\begin{equation}
|b_1|\ge \sqrt{40}\cdot A_0,
\label{bb}
\end{equation}
then
\[
A\le \left(\frac{c_4\cdot H}{A_0}\right)^{1/2}.
\]
For a certain $A_0$ the suitable $H$ is of magnitude $A_0^2$.
A typical
sequence of reduced bounds staring from $A_0=10^{100}$ was the following:
\[
\begin{array}{|c|c|c|c|c|}
\hline
A&10^{100}                        &1.5805\cdot 10^{51}   &6.2833\cdot 10^{26} &1.2528\cdot 10^{15} \\ \hline
H&10^{202}                        &2.4979 \cdot 10^{104} &3.9481\cdot 10^{55} &1.5695\cdot 10^{32} \\   \hline
{\rm new}\; A&1.5805\cdot 10^{51} &6.2833\cdot 10^{26}   &1.2528\cdot 10^{15} &5.5942\cdot 10^8   \\ \hline
\end{array}
\]

\[
\begin{array}{|c|c|c|c|c|}
\hline
A            &5.5942\cdot 10^8    & 3.7382\cdot 10^5    & 8663             &1553            \\ \hline     
H            &3.1295\cdot 10^{19} & 1.3974\cdot 10^{13} & 9.3381\cdot 10^9 &2.4139\cdot 10^8\\ \hline 
{\rm new}\; A & 3.7382\cdot 10^5   & 8663                & 1553            &622             \\ \hline          
\end{array}
\]

\[
\begin{array}{|c|c|c|c|c|}
\hline
A             &622                &  394             &  313             &  280     \\ \hline     
H             & 3.8810\cdot 10^7  & 1.5562\cdot 10^7 &9.8543\cdot 10^6  & 7.8416\cdot 10^6  \\ \hline 
{\rm new}\; A &  394              &  313             & 280              & 264        \\ \hline       
\end{array}
\]

If in a certain step $H$ was not sufficiently large, we replaced it by $10H$.

The reduction procedure was executed with 250 digits accuracy and took only a few seconds.
It has to be performed for each possible values of $j_0$, and the final reduced
bound for $A$ is the maximum of the reduced bounds obtained for $j_0=1,2,3$.

\section{Enumeration, test}

The reduced bound obtained in the previous section gives an upper bound among others for
$|z_{11}|,|z_{12}|$, hence we can enumerate all possible $Z_1$. Further, 
for all possible $\varepsilon,\mu$, equation (\ref{mu}) gives a cubic equation
for $Z_2\in\Z_M$. Testing the roots of this cubic equation in $Z_2$ we can determine
all $Z_2\in\Z_M$ corresponding to $Z_1$.

From (\ref{ZY}) and (\ref{Y}) we can determine 
$Y_1,Y_2$ and then the coordinates $x_{11},x_{12}$ and  $x_{21},x_{22}$ of $X_1,X_2$,
corresponding to $Z_1,Z_2$. Finally, we use (\ref{JG}) to determine $x_{02}$
in the representation (\ref{gamma}) of $\gamma$ (the index of $\gamma$
is independent of $x_{01}$). Substituting the possible tuples
$x_{11},x_{12},x_{21},x_{22}$ into $J(\gamma)$ we obtain a polynomial 
$F(x)=a_9x^9+\ldots +a_1x+a_0$ in $x_{02}$ of degree 9, such that 
\[
|F(x_{02})|=1.
\]
For the roots  $x_{02}$ of absolute value >1 we have
\[
|x_{02}|\le \frac{|a_8|+\ldots |a_1| +  |a_0|+1   }{|a_9|}.
\]
We test the possible integer values of $x_{02}$ and obtain the solutions.
Note that $x_{11},x_{12}$ and  $x_{21},x_{22}$ are usually small values, therefore
the bound for $|x_{02}|$ is also reasonably small.

\section{Example}

We developed and tested our method by taking the trinomial
\[
f(x)=x^6+3x^3+9
\]
with Galois group $C_3\times S_3$ from the paper \cite{hj6} of Harrington and Jones.
These trinomials have several interesting features, which may be the topic
of a separate paper. This polynomial is not monogenic, but the number field $K$
generated  by a root $\alpha$ of it is monogenic.

The quadratic subfield of $K$ is determined by the equation $x^2+3x+9=0$.
It's root is $\beta=(-3+3i\sqrt{3})/2$, therefore $M=\Q(i\sqrt{3})$.
We set $\omega=(1+i\sqrt{3})/2$, then $\beta=3\omega-3$ and $\alpha=\sqrt[3]{\beta}$.
A relative integers basis of $K=\Q(\alpha)$ over $\Q$ is given by
\[
\left(1,\alpha,\frac{\alpha^2(1+\omega)}{3}\right).
\]
We have
\[
\gamma=X_0+X_1\alpha+X_2\frac{\alpha^2(1+\omega)}{3}=Y_0+Y_1\alpha+Y_2\alpha^2,
\]
with
\[
Y_0=X_0,Y_1=X_1,Y_2=X_2\frac{1+\omega}{3}.
\]
Moreover, $\delta=\alpha,k=1,\ell=3$, 
\[
3 \cdot N_{M/\Q}(N_{K/M}(Y_1-\delta Y_2))=\pm 1,
\]
and
\[
N_{M/\Q}(N_{K/M}(Z_1-\delta Z_2))=\pm 3^5
\]
with $Z_1=3Y_1,Z_2=3Y_2$, whence
\[
|Z_1^3-\beta Z_2^3|=|N_{K/M}(Z_1-\delta Z_2))|=3^{5/2}.
\]
Taking $C=10^{100}$ we have to reduce $A$ starting from $24C$. 
The reduction procedure gives a bound 250
for the absolute values of the coordinates $z_{11},z_{12},z_{21},z_{22}$.
In our case $Y_1$ is also integer, hence $z_{11},z_{12}$ are divisible by 3,
which considerable reduces the number of possible pairs $z_{11},z_{12}$.

We used Magma to calculate elements in $\Z_M$ of norm $\pm 3^5$.
We obtained that up to associates the only element is $9-18\omega$.
We set $\varepsilon=\pm 1,\frac{\pm 1\pm i\sqrt{3}}{2}$ and using
\[
Z_1^3-\beta Z_2^3=\varepsilon\mu
\]
we calculated the possible values of $Z_1$, corresponding to $Z_2$.
Finally, we calculated $Y_1,Y_2$, then $X_1,X_2$ and substituted
the coordinates of  $X_1,X_2$ into $J(\gamma)=1$ to determine the 
suitable values of $x_{02}$.  We obtained that up to sign the solutions are:

\[
\begin{array}{|c|c|c|c|c|}
\hline
x_{02}  &  x_{11}  &   x_{12}   &  x_{21}  &   x_{22}   \\ \hline 
-1&0&1&0&-1\\ \hline
1&1&0&1&-1\\ \hline
0&0&0&-1&1\\ \hline
0&0&0&0&-1\\ \hline
0&0&0&1&0\\ \hline
-1&1&-1&1&0\\ \hline
\end{array}
\]
That is, up to sign and translation by $x_{01}$ all generators of power integral bases
of $K$ are given by
\[
\gamma=\omega x_{02}+(x_{11}+ \omega x_{12})\alpha+
(x_{21}+ \omega x_{22})\frac{\alpha^2(1+\omega)}{3},
\]
with the above listed tuples $(x_{02},x_{11},x_{12,}x_{21},x_{22})$.

\end{document}